\documentclass[12pt]{amsart}

\usepackage{citesort}
\usepackage{graphicx}

\newcommand{\0}{\mathbf{0}}
\newcommand{\1}{\mathbf{1}}
\newcommand{\cle}{\preceq}

\newcommand{\bset}{\mathbf{B}}
\newcommand{\cset}{\mathbf{C}}
\newcommand{\pset}{\mathbf{P}}
\newcommand{\rset}{\mathbf{R}}
\newcommand{\maA}{\mathcal{A}}
\newcommand{\CF}{\mathcal{F}}
\newcommand{\rmin}{\mathbf{R}_{\min}}
\newcommand{\rmax}{\mathbf{R}_{\max}}
\newcommand{\Log}{\mathop{\mathrm{Log}}}
\newcommand{\suplim}{\sup\limits}
\newcommand{\sumlim}{\sum\limits}
\newcommand{\maxlim}{\max\limits}

\begin{document}

\title[The Maslov dequantization]
{The Maslov Dequantization, Idempotent and Tropical Mathematics:
a Very Brief Introduction}

\author{G.~L.~Litvinov}
\address{Independent University of Moscow\\ 
Bol'shoi Vlasievskii per., 11\\ 
Moscow 121002, Russia }
\email{islc@dol.ru} 

\date{}

\thanks{This work has been supported by the RFBR grant 05--01--00824
and the Erwin Schr{\"o}\-din\-ger International Institute for
Mathematical Physics (ESI)}

\begin{abstract}
  This paper is a very brief introduction to idempotent mathematics
  and related topics.  It appears as an introductory paper in the
  volume \textit{Idempotent Mathematics and Mathematical Physics}
  (G.~L.~Litvinov and V.~P.~Maslov, eds; AMS Contemporary Mathematics
  Proceedings Series, 2005) \cite{LiMa2005}.
\end{abstract}

\maketitle

This paper is a very brief introduction, without exact theorems and
proofs, to the Maslov dequantization and idempotent and tropical
mathematics.  Our list of references is not complete (not at all).
Additional references can be found, e.g., in the
electronic archive 
\[
  \text{\texttt{http://arXiv.org},}
\]
in \cite{BaCoOlQu92,Bu94,%
Ca79,CoGaQu99,CoQu94,Cu79,Cu95,CuMe80,DeDo98,DeDo2001,%
Fl2004,Gla2002,Gol99,Gol2003,GoMi79,GoMi2001,Gun98a,%
Gun98b,KiRo2004,Ki2001,KoMa97,Ko2001,LiMa95,LiMa98,%
LiMaSh2001,LiSo2001,ZiU81} and in the papers published in this volume.

The author thanks M.\ Akian, Ya.\ I.\ Belopolskaya, P.\ Butkovi{\v c}, G.\ Cohen, S.\ Gaubert,
R.\ A.\ Guninghame-Green, P.\ Del Moral, H.\ Prade, W.\ H.\ Fleming,
J.\ S.\ Golan, M.\ Gondran, I.\ Itenberg, Y.\ Katsov, V.\ N.\ Kolokoltsov,
P.\ Loreti, Yu. I.\ Manin, G.\ Mascari,
W.\ M.\ Mac\-Eneaney, E.\ Pap, M.\ Pedicini, A.\ A.\ Puhalskii, J.-P.\ Quadrat, M.\ A.\ Roytberg, G.\ B.\ Shpiz, I.\ Singer,
and O.\ Viro for useful contacts and references.  Special thanks are due
to V.\ P.\ Maslov for his crucial help and support
and to A.\ N.\ Sobolevski\u\i \ for his great help including, but not
limited to, two pictures presented below.

\section{Some basic ideas}

Idempotent mathematics is based on replacing the usual arithmetic
operations with a new set of basic operations (such as maximum
or minimum), that is on replacing numerical fields by idempotent
semirings and semifields. Typical examples are given by the
so-called max-plus algebra $\rmax$ and the min-plus algebra
$\rmin$. Let $\rset$ be the field of real numbers. Then
${\rmax}={\rset} \cup \{-\infty\}$ with operations $x\oplus y=\max
\{x,y\}$ and $x\odot y=x+y$. Similarly ${\rmin}={\rset}\cup
\{+\infty\}$ with the operations $\oplus=\min$, $\odot=+$. The new
addition $\oplus$ is idempotent, i.e., $x\oplus x=x$ for all
elements $x$.

Many authors (S.\ C.\ Kleene, N.\ N.\ Vorobjev, B.\ A.\ Carre,
R.\ A.\ Cu\-ning\-ha\-me-Gre\-en, K.\ Zimmermann, U.\ Zimmermann, M.\ Gondran, 
F.\ L.\ Baccelli, G.\ Cohen, S.\ Gaubert,
G.\ J.\ Olsder, J.-P.\ Quadrat, and others) used idempotent semirings and
matrices over these semirings for solving some applied problems in
computer science and discrete mathematics, starting from the classical
paper of S.\ C.\ Kleene \cite{Kle56}.  The modern 
\emph{idempotent analysis} (or \emph{idempotent calculus}, or \emph{idempotent
mathematics}) was founded by V.\ P.\ Maslov in the 1980s in Moscow; see,
e.g., \cite{Mas86,Mas87a,Mas87b,MaKo94,MaSa92,MaVo88,KoMa97}.

Idempotent mathematics can be treated as a result of a
dequantization of the traditional mathematics over numerical
fields as the Planck constant $\hbar$ tends to zero taking
imaginary values.  This point of view was presented by 
G.~L.~Litvinov and V.~P.~Maslov \cite{LiMa95,LiMa96,LiMa98}, see also
\cite{LiMaSh2001,LiMaSh2002}.  In other words,
idempotent mathematics is an asymptotic version of the 
traditional
 mathematics over the fields of real and complex numbers.

The basic paradigm is expressed in terms of an {\it idempotent
correspondence principle}.  This principle is closely related to the
well-known correspondence
principle of N.~Bohr in quantum theory.
Actually, there exists a heuristic correspondence between
important, interesting, and useful constructions and results of the
traditional mathematics over fields and analogous constructions and results
over idempotent semirings and semifields (i.e., semirings and semifields
with idempotent addition).

A systematic and consistent application of the idempotent
correspondence principle leads to a variety of results, often
quite unexpected.  As a result, in parallel with the traditional
mathematics over fields, its ``shadow,'' the idempotent
mathematics, appears.  This ``shadow'' stands approximately in the
same relation to the traditional mathematics as does classical physics
to quantum theory, see Fig.~1. 
\begin{figure}[t]
  \centering
  \includegraphics{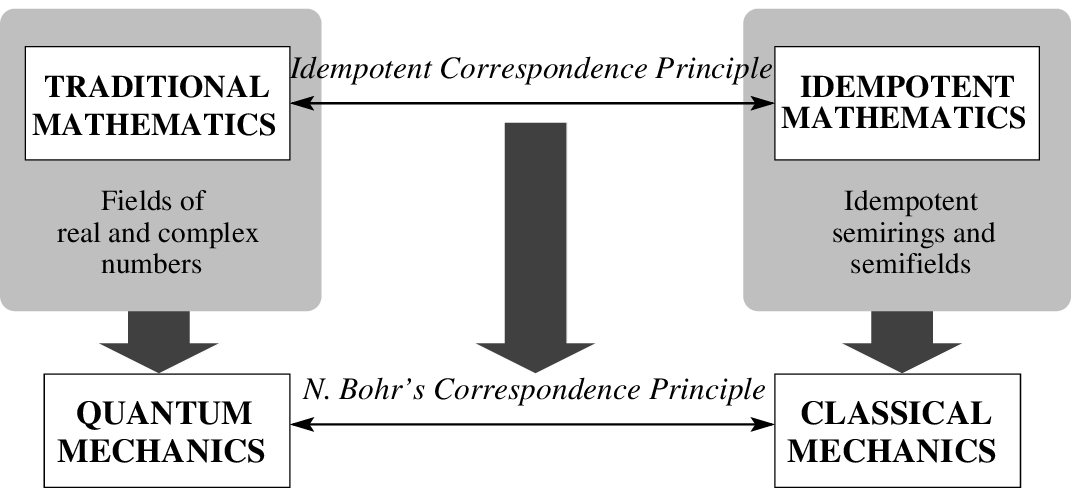}
  \caption{}
\end{figure}
In many respects idempotent mathematics is
simpler than the traditional one.  However the transition from
traditional concepts and results to their idempotent analogs is
often nontrivial.

\section{Semirings, semifields, and dequantization}

Consider a
set $S$ equipped with two algebraic operations: {\it addition} $\oplus$
and {\it multiplication} $\odot$. It is a {\it semiring} if the following
conditions are satisfied:
\begin{itemize}
\item the addition $\oplus$ and the multiplication $\odot$ are
associative;
\item the addition $\oplus$ is commutative;
\item the multiplication $\odot$ is distributive with respect to
the addition $\oplus$: \[x\odot(y\oplus z)=(x\odot y)\oplus(x\odot z)
\quad\text{and}\quad
(x\oplus y)\odot z=(x\odot z)\oplus(y\odot z)\] 
for all $x,y,z\in S$.
\end{itemize}
A {\it unity} of a semiring $S$ is an element $\1\in S$ such that
$\1\odot x=x\odot\1=x$ for all $x\in S$. A {\it zero} of a semiring
$S$ is an element $\0\in S$ such that $\0\neq\1$ and $\0\oplus x=x$,
$\0\odot x=x\odot \0=\0$ for all $x\in S$. A semiring $S$ is called an
{\it idempotent semiring} if $x\oplus x=x$ for all $x\in S$. A
semiring $S$ with neutral elements $\0$ and $\1$ is called a {\it
semifield} if every nonzero element of $S$ is invertible.
Note that dio{\"\i}ds in the sense of \cite{BaCoOlQu92,GoMi79,%
GoMi2001}, quantales in the sense of \cite{Ro90,Ro96},
and inclines in the sense of \cite{KiRo2004} are examples of
idempotent semirings.

Let $\rset$ be the field of real numbers and $\rset_+$ the
semiring of all nonnegative real numbers (with respect to the
usual addition and multiplication). The change of variables $x
\mapsto u = h \ln x$, $h > 0$, defines a map $\Phi_h \colon
\rset_+ \to S = \rset \cup \{-\infty\}$.  Let the addition and
multiplication operations be mapped from $\rset$ to $S$ by
$\Phi_h$, i.e., let $u \oplus_h v = h \ln (\exp (u/h) +
\exp(v/h))$, $u \odot v = u + v$, $\0 = -\infty =$$ \Phi_h(0)$,
$\1 =$$ 0 = \Phi_h(1)$.  It can easily be checked that $u \oplus_h
v \to \max \{u,v\}$ as $h \to 0$ and that $S$ forms a semiring with
respect to addition $u \oplus v = \max \{u,v\}$ and multiplication
$u \odot v = u + v$ with zero $\0 = -\infty$ and unit $\1 = 0$.
Denote this semiring by $\rmax$; it is {\it idempotent}, i.e., $u
\oplus u = u$ for all its elements. The semiring $\rset_{\max}$ is
actually a semifield. The analogy with quantization is obvious;
the parameter $h$ plays the r\^{o}le of the Planck constant, so
$\rset_+$ (or $\rset$) can be viewed as a ``quantum object'' and
$\rmax$ as the result of its ``dequantization.'' A similar
procedure  (for $h<0$) gives the semiring 
$\rmin = \rset \cup \{+\infty\}$ with
the operations $\oplus = \min$, $\odot = +$; in this case $\0 =
+\infty$, $\1$ = 0. The semirings $\rmax$ and $\rset_{\min}$ are
isomorphic.  This passage to $\rmax$ or $\rmin$ is called the
{\it Maslov dequantization}. It is clear that the corresponding
passage from $\cset$ or $\rset$ to $\rset_{\max}$ is generated
by the Maslov dequantization and the map $x\mapsto |x|$.
By misuse of language, {\it we shall also call this passage the
Maslov dequantization}. Connections with physics and
the meaning of imaginary values of the Planck constant are discussed in
\cite{LiMaSh2001,LiMaSh2002}. The idempotent semiring $\rset \cup
\{-\infty\} \cup \{+\infty\}$ with the operations $\oplus = \max$,
$\odot = \min$ can be obtained as a result of a ``second
dequantization'' of $\cset$, $\rset$ or $\rset_+$. Dozens of interesting
examples of nonisomorphic idempotent semirings may be cited as
well as a number of standard methods of deriving new semirings
from these (see, e.g., \cite{CoGaQu2004,Gol99,Gol2003,%
GoMi79,GoMi2001,Gun98a,Gun98b,%
LiMa98,LiMaSh2001}). The so-called
{\it idempotent dequantization} is a generalization of the Maslov
dequantization; this is a passage from fields to idempotent
semifields and semirings in mathematical constructions and
results.

The Maslov dequantization is related to the well-known logarithmic
transformation that was used, e.g., in the classical papers of
E.~Schr{\"o}\-din\-ger (1926) and E.~Hopf (1951). The term `Cole-Hopf
transformation' is also used. The subsequent progress of E.~Hopf's
ideas has culminated in the well-known vanishing viscosity method
and the method of viscosity solutions, see, e.g., \cite{BaCa97,%
CaLi97,FlSo93,MaKo94,Su95} and papers \cite{McC2005} by D.~McCaffrey and 
\cite{Rou2005} by I.~V.~Roublev published in this volume.

\section{Terminology: Tropical semirings and tropical mathematics}

The term `tropical semirings' was introduced in computer
science to denote discrete versions of the max-plus algebra
$\rmax$ or min-plus algebra $\rmin$ 
and their subalgebras; (discrete) semirings of this type
were called tropical semirings by Dominic Perrin in honour of Imre
Simon (who is a Brasilian mathematician and computer
scientist) because of his pioneering activity in this area,
see \cite{Pin98}.

More recently the situation and terminology have changed. For the
most part of modern authors `tropical' means `over $\rset_{\max}$ 
(or $\rset_{\min}$)' and tropical semirings are idempotent
semifields $\rset_{\max}$ and $\rset_{\min}$. The 
terms `max-plus' and `min-plus' are often used in the same
sense. Now the term `tropical mathematics' usually means
`mathematics over $\rset_{\max}$ or $\rset_{\min}$', see,
e.g., \cite{ItKhSh2003,Mi2001,Mi2003,Mi2004,Mi2005,RiStTh2005,SpSt2004}.
Terms
`tropicalization' and `tropification' (see, e.g., \cite{Ki2001})
mean exactly dequantization and quantization in our sense. In
any case, tropical mathematics is a natural and very important
part of idempotent mathematics.

Note that in papers \cite{Vor63,Vor67,Vor70} N.~N.~Vorobjev developed
a version of idempotent linear algebra (with important 
applications, e.g., to mathematical economics) and predicted 
many aspects of the future extended theory.
He used the terms `extremal algebras'
and `extremal mathematics' for idempotent semirings and 
idempotent mathematics.  Unfortunately, 
N.~N.~Vorobjev's papers and ideas were forgotten for a long period, so
his remarkable terminology is not in use any more.

\section{Idempotent algebra and linear algebra}

The first
known paper on idempotent (linear) algebra is due to S.~Kleene
\cite{Kle56}. Systems of linear algebraic equations over an
exotic idempotent semiring of all formal languages over a
fixed finite alphabet are examined in this work; however, 
S.~Kleene's ideas are very general and universal.  Since then, dozens of
authors investigated matrices with coefficients belonging to
an idempotent semiring and the corresponding applications to
discrete mathematics, computer science, computer languages,
linguistic problems, finite automata, optimization problems
on graphs, discrete event systems and Petry nets, stochastic
systems, computer performance evaluation, computational problems
etc. This subject is very well known and well presented in the 
corresponding literature, see, e.g., \cite{BaCoOlQu92,Bu94,%
Bu2005,Ca79,CoGaQu99,CoQu94,Cu79,Cu95,DuSa92,Gla2002,%
Gol99,Gol2003,GoMi79,GoMi2001,Gun98a,Gun98b,%
KiRo2004,KoMa97,Ko2001,LiMa95,LiMa98,LiMa2005,LiMaE2000,%
LiMaRo2000,LiSo2000,LiSo2001,%
MaKo94,Vor63,Vor67,Vor70,ZiU81}. The idempotent linear algebra is
treated in the papers of
P.~Butkovi{\v c} \cite{Bu2005} and E.~Wagneur \cite{Wa2005}  
in the present volume.

Idempotent abstract algebra
is not so well developed yet (on the other hand, from a formal
point of view, the lattice theory and the theory of ordered
groups and semigroups are parts of idempotent algebra). However,
there are many interesting results and applications presented, e.g., 
in \cite{Cu79,Cu95,CuMe80,Ka2004,Ro90,Ro96,Shp2000}.

In particular, an idempotent version of the main theorem of
algebra holds \cite{CuMe80,Shp2000} for radicable idempotent
semifields (a semiring $A$ is {\it radicable} if the equation
$x^n = a$ has a solution $x\in A$ for any $a\in A$ and any
positive integer $n$). It is proved that $\rset_{\max}$ and other
radicable semifields are algebraically closed in a natural sense
\cite{Shp2000}.

\section{Idempotent analysis}

Idempotent analysis was initially
constructed by V.~P.~Maslov and his collaborators and then 
developed by many authors. The subject is presented in the book of
V.~N.~Kolokoltsov and V.~P.~Maslov \cite{KoMa97} (a version of
this book in Russian \cite{MaKo94} was published in 1994).

Let $S$ be an arbitrary semiring with idempotent addition $\oplus$
(which is always assumed to be commutative), multiplication
$\odot$, zero $\0$, and unit $\1$. The set $S$ is supplied with
the {\it standard partial order\/}~$\cle$: by definition, $a \cle
b$ if and only if $a \oplus b = b$. Thus all elements of $S$ are
nonnegative: $\0 \cle$ $a$ for all $a \in S$. Due to the existence
of this order, idempotent analysis is closely related to the lattice
theory, theory of vector lattices, and theory of ordered
spaces. Moreover, this partial order allows to model a number of
basic ``topological'' concepts and results of idempotent analysis at the purely
algebraic level; this line of reasoning was examined
systematically in \cite{LiMaSh98,LiMaSh99,LiMaSh2001,LiMaSh2002,%
LiSh2002} and \cite{CoGaQu2004}.

Calculus deals mainly with functions whose values are numbers. The
idempotent analog of a numerical function is a map $X \to S$,
where $X$ is an arbitrary set and $S$ is an idempotent semiring.
Functions with values in $S$ can be added, multiplied by each
other, and multiplied by elements of $S$ pointwise.

The idempotent analog of a linear functional space is a set of $S$-valued
functions that is closed under addition of functions and multiplication of
functions by elements of $S$, or an $S$-semimodule. Consider, e.g., the
$S$-semimodule $B(X, S)$ of all functions $X \to S$ that are bounded in
the sense of the standard order on $S$.

If $S = \rmax$, then the idempotent analog of integration is defined by the
formula
$$
I(\varphi) = \int_X^{\oplus} \varphi (x)\, dx     = \sup_{x\in X}
\varphi (x),\eqno{(1)}
$$
where $\varphi \in B(X, S)$. Indeed, a Riemann sum of the form
$\sumlim_i \varphi(x_i) \cdot \sigma_i$ corresponds to the expression
$\bigoplus\limits_i \varphi(x_i) \odot \sigma_i = \maxlim_i \{\varphi(x_i)
+ \sigma_i\}$, which tends to the right-hand side of~(1) as $\sigma_i \to
0$. Of course, this is a purely heuristic argument.

Formula~(1) defines the \emph{idempotent} (or \emph{Maslov}) \emph{integral} not only for 
functions taking
values in $\rmax$, but also in the general case when any of bounded
(from above) subsets of~$S$ has the least upper bound.

An \emph{idempotent} (or \emph{Maslov}) \emph{measure} on $X$ is defined by 
$m_{\psi}(Y) = \suplim_{x \in Y}
\psi(x)$, where $\psi \in B(X,S)$. The integral with respect to this
measure is defined by
$$
   I_{\psi}(\varphi)
    = \int^{\oplus}_X \varphi(x)\, dm_{\psi}
    = \int_X^{\oplus} \varphi(x) \odot \psi(x)\, dx
    = \sup_{x\in X} (\varphi (x) \odot \psi(x)).
    \eqno{(2)}
$$

Obviously, if $S = \rmin$, then the standard order $\cle$ is opposite to
the conventional order $\le$, so in this case equation~(2)
assumes the form
$$
   \int^{\oplus}_X \varphi(x)\, dm_{\psi}
    = \int_X^{\oplus} \varphi(x) \odot \psi(x)\, dx
    = \inf_{x\in X} (\varphi (x) \odot \psi(x)),
$$
where $\inf$ is understood in the sense of the conventional order $\le$.

Note that the so-called pseudo-analysis (see a survey paper of
E.~Pap \cite{Pap2005} published in the present volume) is related to a
special part of idempotent analysis; however, this
pseudo-analysis is not a proper part of idempotent 
mathematics in the general case.

\section{Correspondence to stochastics and a duality between 
probability and optimization}

Maslov measures are nonnegative (in the sense of the standard
order) just as probability measures. The analogy between idempotent and
probability measures leads to important relations between
optimization theory and probability theory. By the present time
idempotent analogues of many objects of
stochastic calculus have been constructed and investigated, 
such as  max-plus martingales, max-plus
stochastic differential equations, and others. These results allow, for
example, to transfer powerful stochastic methods
to the optimization theory. This was noticed and examined by
many authors (G.~Salut, P.~Del Moral, M.~Akian, J.-P.~Quadrat,
V.~P.~Maslov, V.~N.~Kolokoltsov, P.~Bernhard, W.~A.~Fleming,
W.~M.~McEneaney, A.~A.~Puhalskii and others), see the survey paper of
W.~A.~Fleming and W.~M.~McEneaney \cite{FlMc2005} published in this
volume and \cite{Ak99,AkQuVi98,Be94,CoQu94,De97,De2004,DeDo98,%
DeDo2001,%
Fl2002,Fl2004,FlMc2000,Gun98a,MaKo94,Puh2001,Qu90,Qu94}.
For relations and applications to large deviations see \cite{Ak99,%
De97,De2004,DeDo98,DeDo2001,Puh2003} and
especially the book of A.~A.~Puhalskii \cite{Puh2001}.

\section{Idempotent functional analysis}

Many other idempotent analogs may be given, in particular, for
basic constructions and theorems of functional analysis.
Idempotent functional
analysis is an abstract version of idempotent analysis. For the
sake of simplicity take $S=\rmax$ and let $X$ be an arbitrary set.
The idempotent integration can be defined by the formula (1), see
above. The functional $I(\varphi)$ is linear over $S$ and its
values correspond to limiting values of the corresponding analogs
of Lebesgue (or Riemann) sums. An idempotent scalar product of
functions $\varphi$ and $\psi$ is defined by the formula
$$
\langle\varphi,\psi\rangle = \int^{\oplus}_X \varphi(x)\odot\psi(x)\, dx =
\sup_{x\in X}(\varphi(x)\odot\psi(x)).
$$
So it is natural to construct idempotent analogs of integral
operators in the form
$$
K : \varphi(y) \mapsto (K\varphi)(x) = \int^{\oplus}_Y K(x,y)\odot
\varphi(y)\, dy = \sup_{y\in Y}\{K(x,y)+\varphi(y)\},\eqno(3)
$$
where $\varphi(y)$ is an element of a space of functions defined
on a set $Y$, and $K(x,y)$ is an $S$-valued function on $X\times Y$.
Of course, expressions of this type are standard in optimization
problems.\medskip

Recall that the definitions and constructions described above can
be extended to the case of idempotent semirings which are
conditionally complete in the sense of the standard order. Using
the Maslov integration, one can construct various function spaces
as well as idempotent versions of the theory of generalized
functions (distributions). For some concrete idempotent function
spaces it was proved that every `good' linear operator (in the
idempotent sense) can be presented in the form (3); this is an
idempotent version of the kernel theorem of L.~Schwartz; results
of this type were proved by V.~N.~Kolokoltsov, P.~S.~Dudnikov and
S.~N.~Samborski, I.~Singer, M.~A.~Shubin and others, see, e.g., 
\cite{DuSa92,KoMa97,MaKo94,MaSa92,Si2003}. So every
`good' linear functional can be presented in the form
$\varphi\mapsto\langle\varphi,\psi\rangle$, where
$\langle,\rangle$ is an idempotent scalar product.\medskip

In the framework of idempotent functional analysis results of this
type can be proved in a very general situation. In \cite{LiMaSh98,%
LiMaSh99,LiMaSh2001,LiMaSh2002,LiSh2002} an algebraic version 
of the idempotent functional
analysis is developed; this means that basic (topological) notions
and results are simulated in purely algebraic terms. The treatment
covers the subject  from basic concepts and results (e.g.,
idempotent analogs of the well-known theorems of Hahn-Banach,
Riesz, and Riesz-Fisher) to idempotent analogs of 
A.~Grothendieck's concepts and results on topological tensor
products, nuclear spaces and operators. An abstract version of the
kernel theorem is formulated. Note that the passage from the usual
theory to idempotent functional analysis may be very nontrivial;
for example, there are many non-isomorphic idempotent Hilbert
spaces. Important results on idempotent functional analysis
(duality and separation theorems) are recently published by 
G.~Cohen, S.~Gaubert, and J.-P.~Quadrat \cite{CoGaQu2004};
see also a finite dimensional version of the separation
theorem in \cite{ZiK77}.
Three papers on this subject by M.~Akian, S.~Gaubert, V.~Kolokoltsov,
G.~Cohen, J.-P.~Quadrat, I.~Singer, and C.~Walsh
\cite{AkGaKo2005,AkGaWa2005,CoGaQuSi2005} are published in this
volume.

There is an ``idempotent'' version
of the theory of linear representations of groups and semigroups
and the abstract
harmonic analysis, see, e.g., \cite{LiMaSh2002}. 
In the framework of this theory
the well-known Legendre transform can be treated as an idempotent
version of the traditional Fourier transform (this observation is
due to V.~P.~Maslov).\medskip

\section{The superposition principle and linear problems}

Basic
equations of quantum theory are linear; this is the superposition
principle in quantum mechanics. The Hamilton--Jacobi equation,
the basic equation of classical mechanics, is nonlinear in the
conventional sense. However, it is linear over the semirings
$\rmax$ and $\rmin$. Similarly, different versions of the Bellman
equation, the basic equation of optimization theory, are linear
over suitable idempotent semirings; this is V.~P.~Maslov's idempotent
superposition principle, see \cite{Mas86,Mas87a,Mas87b,MaKo94,%
MaSa92}. For instance, the
finite-dimensional stationary Bellman equation can be written in
the form $X = H \odot X \oplus F$, where $X$, $H$, $F$ are
matrices with coefficients in an idempotent semiring $S$ and the
unknown matrix $X$ is determined by $H$ and $F$. In particular,
standard problems of dynamic programming and the well-known
shortest path problem correspond to the cases $S = \rmax$ and $S
=\rmin$, respectively. It is known that principal optimization
algorithms for finite graphs correspond to standard methods for
solving systems of linear equations of this type (i.e., over
semirings). Specifically, Bellman's shortest path algorithm
corresponds to a version of Jacobi's algorithm, Ford's algorithm
corresponds to the Gauss--Seidel iterative scheme, etc.

The linearity of the Hamilton--Jacobi equation over
$\rmin$ and $\rmax$, which is the result of
the Maslov dequantization of the Schr{\"o}\-din\-ger equation, 
is closely related to the (conventional)
linearity of the Schr{\"o}\-din\-ger equation and can be deduced from
this linearity. Thus, it is possible to borrow standard ideas and
methods of linear analysis and apply them to a new area.

The action functional $S=S(x(t))$ can be considered as a function
taking the set of curves (trajectories) to the set of real numbers
which can be treated as elements of  $\rmin$. In this case the
minimum of the action functional can be viewed as the Maslov
integral of this function over the set of trajectories or an
idempotent analog of the Euclidean version of the Feynman
path integral. The minimum of the
action functional corresponds to the maximum of $e^{-S}$, i.e.
idempotent integral $\int^{\oplus}_{\{paths\}} e^{-S(x(t))}
D\{x(t)\}$. Thus the least action principle can be considered as
an idempotent version of the well-known Feynman approach to
quantum mechanics.  The representation of a solution to the
Schr{\"o}\-din\-ger equation in terms of the Feynman integral 
corresponds to the Lax--Ole\u{\i}nik solution formula for the
Hamilton--Jacobi equation.

The idempotent superposition principle indicates that there exist important
nonlinear (in the traditional sense) problems that are linear over
idempotent semirings.  The linear idempotent functional analysis is a natural
tool for investigation of those nonlinear infinite-dimensional problems that
possess this property.

\section{Dequantization of geometry}

An idempotent version of real
algebraic geometry was discovered in the report of O.~Viro for the
Barcelona Congress \cite{Vir2000}. Starting from the idempotent 
correspondence principle O.~Viro constructed a piecewise-linear
geometry of polyhedra of a special kind in finite dimensional 
Euclidean spaces as a result of the
Maslov dequantization of real algebraic geometry. He indicated
important applications in real algebraic geometry (e.g., 
 in the framework of Hilbert's 16th problem for
constructing real algebraic varieties with prescribed properties
and parameters) and relations to complex algebraic geometry and
amoebas in the sense of I.~M.~Gelfand, M.~M.~Kapranov, and
A.~V.~Zelevinsky (see their book \cite{GeKaZe94} and \cite{Vir2002}).
Then complex algebraic geometry was dequantized by 
G.~Mikhalkin and the result turned out
to be the same; this new `idempotent' (or
asymptotic) geometry is now often called the {\it tropical algebraic geometry},
see, e.g., \cite{EiKaLi2004,ItKhSh2003,Mi2001,Mi2003,Mi2004,%
Mi2005,RiStTh2005,Shu2002,SpSt2004,St2002}.

There is a natural relation between the Maslov dequantization
and amoebas. 
Suppose $({\cset}^*)^n$ is a complex torus, where
${\cset}^* = {\cset}\backslash \{0\}$ is the group
of nonzero complex numbers under multiplication.  For
 $z = (z_1, \dots, z_n)\in 
(\cset^*)^n$ and a positive real number $h$ denote by
$\Log_h(z) = h\log(|z|)$ the element
\[(h\log |z_1|, h\log |z_2|, \dots,
h\log|z_n|) \in \rset^n.\]  
Suppose $V\subset (\cset^*)^n$ is a
complex algebraic variety; denote by $\maA_h(V)$ the set
$\Log_h(V)$. If $h=1$, then the set $\maA(V) = \maA_1(V)$ is
called the {\it amoeba} of $V$ in the sense of \cite{GeKaZe94},
see also \cite{Vir2002};
the amoeba $\maA(V)$ is a closed subset of $\rset^n$ with a
non-empty complement. 
Note that this construction depends on our coordinate
system.

For the sake of simplicity suppose $V$ is a hypersurface 
in~$(\cset^*)^n$ defined by a polynomial~$f$; then there 
is a deformation $h\mapsto f_h$ of this polynomial generated by the
Maslov dequantization and $f_h = f$ for $h = 1$.
Let $V_h\subset ({\cset}^*)^n$ be the zero set of $f_h$ and
set $\maA_h (V_h) = {\Log}_h (V_h)$. Then
 there exists a tropical variety
$\mathit{Tro}(V)$ such that the subsets $\maA_h(V_h)\subset \rset^n$
tend to $\mathit{Tro}(V)$ in the Hausdorff metric as $h\to 0$, 
see \cite{Mi2001,Ru2001}.
The tropical variety $\mathit{Tro}(V)$ is a result
of a deformation of the amoeba $\maA(V)$ and the Maslov
dequantization of the variety $V$. The set $\mathit{Tro}(V)$ is called the
{\it skeleton} of $\maA(V)$. 

\begin{figure}
  \centering
  \includegraphics{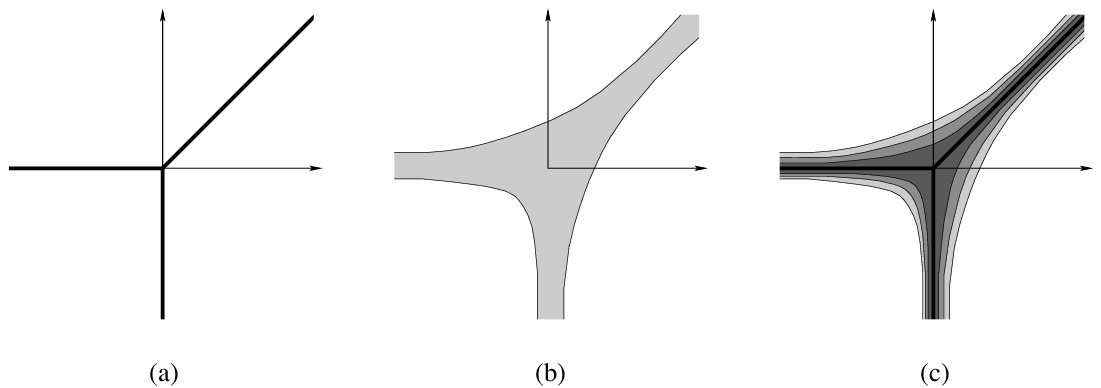}
  \caption{}
\end{figure}

{\bf Example \cite{Mi2001}.}
 For the line $V = \{\, (x, y)\in ({\cset}^*)^2 \mid
x + y + 1 = 0\,\}$ the piecewise-linear graph $\mathit{Tro}(V)$ is a tropical
line, see Fig.~2(a). The amoeba $\maA(V)$ is represented in Fig.~2(b), while
Fig.~2(c) demonstrates the corresponding deformation
of the amoeba.

There is an interesting paper \cite{PaTs2005} of M.~Passare and A.~Tsikh on
amoebas of algebraic and analytic varieties in the present
volume.

In the important paper \cite{Kap2000} (see also 
\cite{EiKaLi2004,Mi2001,Mi2004,RiStTh2005}) tropical
varieties appeared as amoebas over non-Archimedean fields. In
2000 M.~Kontsevich noted that it is possible to use non-Arhimedean
amoebas in enumerative geometry, see \cite[section
2.4, remark 4]{Mi2001}. In fact methods of tropical geometry lead to
remarkable applications to the algebraic enumerative geometry,
Gromov-Witten and Welschinger invariants, see \cite{ItKhSh2003,%
ItKhSh2004,Mi2001,Mi2003,Mi2004,Mi2005}. In particular,
G.~Mikhalkin presented and proved in \cite{Mi2003,Mi2005} a
formula enumerating curves of arbitrary genus in toric
surfaces.

Note that tropical geometry is closely related to the
well-known program of M.~Kontsevich and Y.~Soibelman, see,
e.g., \cite{KonSoi2001,KonSoi2004}. There is an introductory
paper \cite{RiStTh2005} on tropical algebraic geometry in this
volume. The paper of G.~L.~Litvinov and G.~B.~Shpiz \cite{LiSh2005}
(which is also published in the present volume) is also related to
the subject.

However on the whole only first steps in idempotent/tropical geometry
have been made and the problem of systematic construction of idempotent
versions of algebraic and analytic geometries is still open.

\section{The correspondence principle for algorithms and their computer
implementations}

There are many important applied algorithms of idempotent
mathematics, see, e.g., \cite{BaCoOlQu92,Bu2005,Ca79,Cu95,CuMe80,FiRo93,%
FlMc2000,FlMc2005,GoMi79,GoMi2001,Gun98a,ItKhSh2003,KiRo2004,Ki2001,%
Ko2001,KoMa97,LiMa95,LiMa98,LiMaE2000,LiMaRo2000,LiSo2000,LiSo2001,%
LoPe2005,Mi2003,Mi2005,RiStTh2005,Roy92,St2002,Vor63,Vor67,Vor70,%
WoOl2005,ZiK2005,ZiU81}. The idempotent correspondence principle is valid
for algorithms as well as for their software and hardware implementations
\cite{LiMa95,LiMa96,LiMa98,LiMaE2000,LiMaRo2000}. 
In particular, according to the superposition principle, analogs
of linear algebra algorithms are especially important. It is well-known
that algorithms of linear algebra are convenient for parallel computations;
so their idempotent analogs accept a parallelization. This is a regular way
to use parallel computations for many problems including basic optimization
problems. It is convenient to use universal algorithms which do not
depend on a concrete semiring and its concrete computer model. Software
implementations for universal semiring algorithms are based on
object-oriented and generic programming; program modules can deal with
abstract (and variable) operations and data types,
see \cite{LiMa95,LiMa98,LiMaE2000,LiMaRo2000}. The paper
\cite{LoPe2005} of P.~Loreti and M.~Pedicini on the subject
is published in the present volume.

The most important and standard algorithms have many hardware
implementations in the form
of technical devices or special processors. These devices often can be used as
prototypes for new hardware units generated by substitution of the
usual arithmetic operations for its semiring analogs, see 
\cite{LiMa95,LiMa98,LiMaRo2000}.
Good and efficient technical ideas and decisions can be
transposed from prototypes into new hardware units. Thus the correspondence
principle generates a regular heuristic method for hardware design.

\section{Idempotent interval analysis}

An idempotent version of the
traditional interval analysis is presented in \cite{LiSo2000,LiSo2001}.
Let $S$ be an idempotent semiring equipped with the
standard partial order. A {\it closed interval} in $S$ is a subset
of the form ${\bf x} = [{\bf \underline x}, {\bf \bar  x}] = \{
x\in S\mid {\bf \underline x}\cle x
\cle{\bf \bar x}\}$, where the elements $\bf \underline x\cle \bf
\bar x$ are called {\it lower} and {\it upper bounds} of the
interval $\bf x$. A {\it weak interval extension} $I(S)$ of the
semiring $S$ is the set of all closed intervals in $S$ endowed
with operations $\oplus$ and $\odot$ defined as $\bf x\oplus\bf y
= [\bf \underline x \oplus\bf \underline y, \bf \bar
x\oplus\bf\bar y]$, $\bf x\odot\bf y = [\bf \underline x\odot\bf
\underline y, \bf \bar x\odot\bf \bar y]$; the set $I(S)$ is a new
idempotent semiring with respect to these operations.  It is
proved that basic interval problems of idempotent linear algebra
are polynomial, whereas in the traditional interval analysis
problems of this kind are generally NP-hard. Exact interval
solutions for the discrete stationary Bellman equation (see the
matrix equation discussed in section 8 above) and for the 
corresponding optimization problems are constructed and examined by
G.~L.~Litvinov and A.~Sobolevski{\u\i}
in \cite{LiSo2000,LiSo2001}. Similar results are presented by
K.~Cechl{\'a}rov{\'a} and R.~A.~Cu\-ning\-ha\-me-Gre\-en
in~\cite{CeCu2002}.

\section{Relations to the KAM theory and optimal transport}

The subject of the
Kolmogorov--Arnold--Moser (KAM) theory may be formulated as the
study of invariant subsets in phase spaces of nonintegrable
Hamiltonian dynamical systems where the dynamics displays the same
degree of regularity as that of integrable systems (quasiperiodic
behaviour). Recently, a considerable progress was made via a
variational approach, where the dynamics is specified by the
Lagrangian rather  than Hamiltonian function. The corresponding
theory was initiated by S.~Aubry and J.~N.~Mather and recently
dubbed {\it weak KAM theory} by A.~Fathi (see his book ``Weak KAM
Theorems in Lagrangian Dynamics,'' Cambridge Univ.\ Press, 2004; see
also \cite{KKS2005a,KKS2005b,So99a,So99b}). 
Minimization of a certain
functional along trajectories of moving particles is a central
feature of another subject, the optimal transport theory, which
also has undergone a rapid recent development. This theory dates
back to G.~Monge's work on cuts and fills (1781). A modern version
of the theory is known now as the so-called {\it
Monge--Amp\`ere--Kantorovich (MAK)  optimal transport theory}
(after the work of L.~V.~Kantorovich ``On the translocation of
masses'' in C.R. (Doklady) Acad. Sci. USSR, v. 321, 1942,
p.~199--201). There is a similarity between the two theories and
there are relations to problems of the idempotent functional
analysis (e.g., the problem of eigenfunctions for ``idempotent''
integral operators, see \cite{So99a}). Applications of optimal
transport to data processing in cosmology are presented in
\cite{BFHLMMS2003,FMMS2002}.

\section{Relations to logic, fuzzy sets, and possibility theory}

Let $S$ be an idempotent semiring with neutral elements $\0$ and
$\1$ (recall that $\0\neq \1$, see section 2 above). Then the Boolean
algebra $\bset = \{\0, \1\}$ is a natural idempotent subsemiring
of $S$. Thus $S$ can be treated as a generalized (extended) logic
with logical operations $\oplus$ (disjuction) and $\odot$ 
(conjunction). Ideas of this kind are discussed in many books and
papers with respect to generalized versions of logic and especially
 quantum logic, see, e.g., \cite{Gol99,KlPa2004,Ro90,Ro96}. In 
the present volume there is a paper of A.~Di~Nola and B.~Gerla
\cite{DiGe2005} related to these ideas.

Let $\Omega$ be the so-called universe consisting of ``elementary events.''
Denote by ${\CF}(S)$ the set of functions defined on $\Omega$ and taking
their values in $S$; then ${\CF}(S)$ is an idempotent semiring with respect
to the pointwise addition and multiplication of functions. We shall
say that elements of ${\CF}(S)$ are {\it generalized fuzzy sets}, 
see \cite{Gol99,Li2004}. We have the 
well-known classical definition of fuzzy sets (L.~A.~Zadeh \cite{Za65})
if $S = {\pset}$, where $\pset$ is the segment $[0,1]$ with the semiring operations
$\oplus = \max$ and $\odot = \min$. Of course,
functions from ${\CF}(\pset)$ taking their values in the Boolean algebra
$\bset = \{0, 1\}\subset {\pset}$ correspond to traditional sets from $\Omega$
and semiring operations correspond to standard operations for sets. In the
general case 
functions from ${\CF}(S)$ taking their values in $\bset = \{{\0}, 
{\1}\}\subset S$ can be treated as traditional subsets in $\Omega$. If $S$
is a lattice (i.e. $x\odot y = \inf \{x, y\}$ and $x\oplus y = \sup 
\{x, y\}$), then generalized fuzzy sets coincide with $L$-fuzzy sets in the
sense of J.~A.~Goguen \cite{Go67}. The set $I(S)$ of intervals is an idempotent
semiring (see section 11), so elements of ${\CF}(I(S))$ can be treated as
interval (generalized) fuzzy~sets.

It is well known that the classical theory of fuzzy sets is a basis
for the theory of possibility \cite{Za78,DuPrSa2001}. Of course, it is possible
to develop a similar generalized theory of possibility starting from
generalized fuzzy sets, see, e.g., \cite{DuPrSa2001,KlPa2004,Li2004}.
Generalized theories can be 
noncommutative; they seem to be more qualitative and less quantitative
with respect to the classical theories presented in \cite{Za65,Za78}. 
We see that
idempotent analysis and the theories of (generalized) fuzzy sets and 
possibility have the same objects, i.e. functions taking their values in
semirings. However, basic problems and methods could be different
for these theories (like for the measure theory and the probability
theory).

\section{Relations to other areas and miscellaneous applications}

Many relations and applications of idempotent mathematics to different
theoretical and applied areas of mathematical sciences are discussed
above. Of course, optimization and optimal control problems form a very natural 
field for applications of ideas and methods of idempotent mathematics. There
is a very good survey paper \cite{Ko2001} of V.~N.~Kolokoltsov on the
subject, see also \cite{BaCoOlQu92,Bu2005,Ca79,CoGaQu99,CoQu94,Cu79,Cu95,%
CuMe80,De97,DeDo98,DeDo2001,Fl2002,Fl2004,FlMc2000,FlMc2005,%
GoMi79,GoMi2001,Gun98a,LiMa95,LiMa98,LiSo2000,LiSo2001,LoQuMa2005,%
Mas86,%
Mas87a,Mas87b,MaKo94,MaSa92,MaVo88,Qu94,Vor63,Vor67,Vor70,%
WoOl2005,ZiK2005,ZiU81}.

There are many applications to differential equations and
stochastic differential equations, see, e.g., 
\cite{Fl2002,Fl2004,FlMc2000,FlMc2005,Gun98a,Ko96,Ko2000,Ko2001a,%
KoMa97,Mas86,Mas87a,Mas87b,MaKo94,MaSa92,Pap2005,So99a,So99b}.

Applications to the game theory are discussed, e.g., in 
\cite{KoMal97,KoMa97,MaKo94}. There are interesting applications
in biology (bioinformatics), see, e.g., \cite{FiRo93,PaSt2004,Roy92}. Applications
and relations to mathematical morphology are examined the paper 
\cite{DeDo2001} of P.~Del Moral and M.~Doisy and especially in an extended
preprint version of this article. There are many relations and
applications to physics (quantum and classical physics, 
statistical physics, cosmology etc.) see, e.g., \cite{Ko2000,%
KoMa97,LiMaSh2001,LiMaSh2002,Nu91,Qu97,ChDu87}, section 12 above
and the paper of P.~Lotito, J.-P.~Quadrat, E.~Mancinelli
\cite{LoQuMa2005} published in this volume.

There are important relations and applications to purely mathematical areas.
The so-called tropical combinatorics is discussed in a large survey paper
\cite{Ki2001} of A.~N.~Kirillov, see also \cite{Bu2005,ZiU81}. 
Tropical mathematics is closely
related to the very attractive and popular theory of cluster
algebras founded by S.~Fomin and A.~Zelevinsky, see their survey
paper \cite{FoZe2004}. In both cases there are relations with the
traditional theory of representations of Lie groups and related topics.
There are important relations with convex analysis and discrete
convex analysis, see, e.g., \cite{AkGa2003,CoGaQuSi2005,Cu95,%
DeSt2003,LiSh2005,MaTi2003,MaSa92,Si97,ZiK77,ZiK79a,ZiK79b} 
and the paper of G.~Cohen, S.~Gaubert, 
J.-P.~Quadrat, and J.~Singer published in the present
volume.

Many authors examine, explicitly or not, relations and applications of idempotent
mathematics to mathematical economics starting from the 
classical papers of N.~N.~Vorobjev \cite{Vor63,Vor67,Vor70}, see,
e.g., \cite{DaKoMu2001,KoMal97,MaSa92,ZiK76,ZiU81}.

\end{document}